\documentclass{amsart}
\usepackage{amsmath, amssymb}

\newtheorem{thm}{Theorem}[section]
\newtheorem{cor}[thm]{Corollary}
\newtheorem{lem}[thm]{Lemma}
\newtheorem{prop}[thm]{Proposition}
\theoremstyle{definition}

\theoremstyle{remark}
\newtheorem{rem}[thm]{Remark}
\newcommand{\Z}{{\mathbb Z}}
\newcommand{\R}{{\mathbb R}}

\newcommand{\C}{{\mathbb C}}
\newcommand{\CP}{{\mathbb{CP}}}
\newcommand{\CS}{Y}
\newcommand{\CSt}{\widetilde{\CS}}
\newcommand{\Mcg}{\operatorname{Mod}}
\newcommand{\PSL}{\mathrm{PSL}}
\newcommand{\PU}{\mathrm{PU}}

\newcommand{\F}{\mathbb{F}}

\newcommand{\ord}{\operatorname{ord}}
\numberwithin{equation}{section}

\begin{document}
\title{Braid monodromy and  the Cartwright--Steger surface}
\author{Andr\'{a}s I. Stipsicz}
\address{HUN-REN R\'enyi Institute of Mathematics\\
H-1053 Budapest\\ 
Re\'altanoda utca 13--15, Hungary}
\email{stipsicz.andras@renyi.hu}

\author{Zolt\'an Szab\'o}
\address{Department of Mathematics\\
Princeton University\\
 Princeton, NJ, 08544}
\email{szabo@math.princeton.edu}

\begin{abstract}
  Building on the Livné-like construction of Cartwright, Koziarz and Yeung, we construct a closed
  smooth \(4\)-manifold with a map to \(\mathbb{CP}^1\)
  whose generic fiber has genus \(109\) and whose four
  exceptional fibers are determined by explicit braid and permutation data.
  A pencil of four-pointed spheres realizes the
  same global braid monodromy and makes the covering holomorphic. We show that the contraction
  of its thirty-six exceptional curves is the Cartwright–Steger surface, hence our smooth construction
  gives a topological description of this complex surface.  
\end{abstract}
\maketitle

\section{Introduction}
The Cartwright--Steger surface is a smooth compact quotient of the complex
two-ball with Euler number $3$ and first Betti number $b_1=2$ \cite{CS}.  Three is the
smallest Euler number of a smooth surface of
general type.  Among the surfaces attaining it, the fake projective
planes and the Cartwright--Steger surface, the latter is the unique
example with nonzero first Betti number and the unique one admitting
a real structure \cite{PY,CS,KK,BY}.
Cartwright--Koziarz--Yeung constructed on a thirty-six-fold blow-up of
this surface a Livn\'e-like fibration whose generic fiber is a
$288$-fold branched cover of a four-pointed sphere and has genus
$109$ \cite{CKY}.

In Section~\ref{sec:A} we construct a genus-$109$ fibration from an
explicit braid representation and a finite permutation tuple.
 Section~\ref{sec:B} realizes
the same global braid representation by a rational pencil, constructs the
holomorphic cover, and identifies its blow-down with the
Cartwright--Steger surface by a canonical-divisor calculation.

In more concrete terms,
using braid monodromy techniques we construct
a four-manifold $Z^\circ$ together with a map
\[
Z^\circ \to \CP^1_\circ
\]
where $\CP^1_\circ$ denotes the complement of the four points
$\{ 0,1,\omega , \omega ^2\}\subset \CP ^1$ with $\omega = e^{2\pi i/3}$
and the following
properties:

\begin{thm}\label{thm:A}
  The bundle $Z^\circ\to\CP^1_\circ$
  extends over the four missing values to a closed smooth
  $4$-manifold $Z$ with a smooth map $Z\to\CP^1$, with
  the following properties.
  \begin{enumerate}
  \item Over $\CP^1_\circ$ it is a fiber bundle with
    fiber a closed orientable surface of genus $109$.
  \item Over each of $1,\omega,\omega^2$ the fiber has multiplicity $2$, with
    smooth reduction of genus $55$.
  \item Over $0$ the fiber has four components whose normalizations have
    genera $4,4,10,10$.  Their branches
    meet at three points, at each of which eight branches cross pairwise
    transversally.  The Euler characteristic of this fiber is $-69$.
  \item The Euler characteristic of the closed four-manifold $Z$ is
    $\chi(Z)=39$.
  \item There are thirty-six disjoint smoothly embedded spheres
    $\Sigma_1,\dots,\Sigma_{36}\subset Z$, each of self-intersection $-1$.
       Blowing the spheres down gives a closed smooth $4$-manifold $Z_0$
    with $\chi(Z_0)=3$.
  \end{enumerate}
\end{thm}

In Section~\ref{sec:B} we provide an alternative construction which
makes the fibration of Theorem \ref{thm:A} holomorphic and identifies
its blow-down with the Cartwright--Steger surface.  The input is an explicit
rational family of four-pointed spheres:  In $\CP^2$ with coordinates
$[x_0:x_1:x_2]$ let
\begin{equation}\label{eq:pencildata}
   \mathcal C\colon x_0^3=x_1^2x_2,
   \quad
   V\colon x_1=x_2,
   \quad
   q_0=[0:-\tfrac12:1],
   \quad
   L_z\colon 3x_0=z(2x_1+x_2),
\end{equation}
so that $\mathcal C$ is a cuspidal cubic,  $\{ L_z\} _{z\in\CP^1}$
is the pencil of lines through $q_0$, and $V$ is the line through the
three points at which a line of the pencil is tangent to $\mathcal C$.
On each $L_z$ mark the point $L_z\cap V$, which we distinguish,
together with the three points of $L_z\cap\mathcal C$.

\begin{thm}\label{thm:braidmatch}
  The four marked points of $L_z$ are distinct precisely for
  $z\in\CP^1_\circ$.  At $z=0$, two cubic intersection points meet
  at the cusp with local behavior $\pm z^{3/2}$.  At each of
  $1,\omega,\omega^2$, a tangency occurs on $V$.
    With suitable fiber marking and meridian paths, the
  global braid representation is exactly the representation
  used in the proof of Theorem~\ref{thm:A}.
\end{thm}

The 288-fold branched cover appearing in the smooth construction now
can be applied in the holomorphic setting (with appropriately chosen
branch curves), providing a smooth complex surface with the following
properties:

\begin{thm}\label{thm:holreal}
  \begin{enumerate}
  \item Taking in each fiber of the pencil of Equation~\eqref{eq:pencildata}
    a suitable $288$-fold branched
    cover, and completing over
    $0,1,\omega,\omega^2$, produces a smooth compact complex surface $W$ with a
    holomorphic map $W\to\CP^1$.
     \item $W\to\CP^1$ is fiber-preservingly diffeomorphic to the fibration
    $Z\to\CP^1$ of Theorem \ref{thm:A}.
  \item Complex conjugation preserves the covering monodromy class and
    therefore lifts to an antiholomorphic involution $\tilde\kappa$ of $W$ covering
    $z\mapsto\bar z$ on the base.  It preserves the collection of thirty-six disjoint
    $(-1)$-spheres of Theorem \ref{thm:A}(5) and hence descends to a real
    structure $c$ on their blow-down $W_0$.  Exactly four of the thirty-six
    contracted spheres are invariant, each with nonempty real locus, and
    \[
      \chi\bigl(\operatorname{Fix}(c)\bigr)=-1,
      \qquad
      \chi\bigl(\operatorname{Fix}(\tilde\kappa)\bigr)=-5 .
    \]
  \end{enumerate}
\end{thm}

The real structure on $W$ and $W_0$ above then allows us
to identify the blow-down surface $W_0$ (and hence $Z_0$)
with the Cartwright--Steger surface. It also shows that the real
structure $c$ on $W_0$ is holomorphically conjugate to the one constructed by
Borisov--Yeung \cite{BY}.  More precisely:
\begin{thm}\label{thm:CSid}
  The blow-down $W_0$ is biholomorphic to the Cartwright--Steger
  surface $\CS$.  A biholomorphism may be chosen to identify its
  real structure $c$ with the Borisov--Yeung real structure $c_{BY}$ on
  the Cartwright-Steger surface $Y$.
\end{thm}

\subsection{Acknowledgments}
The first author was supported by ERC Advanced Grant KnotSurf4d.  The
second author was partially supported by the Simons Grant \emph{New
Structures in Low-Dimensional Topology}.  We thank \.Inan\c{c} Baykur
for helpful conversations.  LLM's were used in the preparation of this
paper; in particular, we used Python programs written by AI tools to
check identities regarding representations.  LLM's were also used to
help in typesetting and language editing.  All mathematical content
has been checked by the authors, who take full responsibility for it.

\section{The smooth construction}\label{sec:A}

\subsection{The construction of the bundle with generic fiber the
    four-punctured sphere}
  Let $S$ be the 2-dimensional sphere with one distinguished marked
  point $B$ and three further marked points $A_0,A_1,A_2$.  First we
  consider an $S$-bundle (each fiber coming with four marked points)
  over $\CP^1_\circ =\CP ^1\setminus \{ 0,1,\omega , \omega ^2\}$.

Write $\Mcg(S;B,\{A_i\})$ for the group of
orientation-preserving diffeomorphisms of $S$ that fix $B$ and
preserve $\{A_0,A_1,A_2\}$ setwise, modulo isotopy through such
diffeomorphisms.  Removing $B$ from $S$ identifies
$$\Mcg\bigl(S;B,\{A_i\}\bigr)\;\cong\;B_3/Z(B_3)\;\cong\;\PSL_2(\Z),$$
the mapping class group of a disc with three marked points modulo
rotation of the boundary.
Write $\bar\sigma_1,\bar\sigma_2$ for the half-twists,
$\bar\Delta=\bar\sigma_1\bar\sigma_2\bar\sigma_1$ and
$\varrho=\bar\sigma_1\bar\sigma_2$, so that $\bar\Delta^2=\varrho^3=1$; under
$\bar\sigma_1\mapsto\left(\begin{smallmatrix}1&1\\0&1\end{smallmatrix}\right)$,
$\bar\sigma_2\mapsto\left(\begin{smallmatrix}1&0\\-1&1\end{smallmatrix}\right)$
one has
$\bar\Delta\mapsto\left(\begin{smallmatrix}0&1\\-1&0\end{smallmatrix}\right)$ and
$\varrho\mapsto\left(\begin{smallmatrix}0&1\\-1&1\end{smallmatrix}\right)$.

\begin{lem}\label{lem:relation}
Consider the elements
  \[
    \beta_0=\bar\sigma_1^{\,3},\quad \beta_1=\bar\Delta,\quad
    \beta_\omega=\varrho\bar\Delta\varrho^{-1},\quad
    \beta_{\omega^2}=\varrho^2\bar\Delta\varrho^{-2}.
  \]
  Then $\beta_0\beta_1\beta_\omega\beta_{\omega^2}=1$.
\end{lem}
\begin{proof}
  Since $\varrho^3=1$,
  $\beta_1\beta_\omega\beta_{\omega^2}=(\bar\Delta\varrho)^3$, and
  $\bar\Delta\varrho=\left(\begin{smallmatrix}-1&1\\0&-1\end{smallmatrix}\right)
    =\bar\sigma_1^{-1}$ in $\PSL_2(\Z)$, so the product is
    $\left(\begin{smallmatrix}1&-3\\0&1\end{smallmatrix}\right)=\bar\sigma_1^{-3}$.
\end{proof}

Since $\pi_1$ of the four-punctured sphere is the free group
generated by three of the
four
loops, with the product of all four trivial, Lemma \ref{lem:relation} defines a
homomorphism
$$\beta\colon\pi_1\bigl(\CP^1\setminus\{0,1,\omega,\omega^2\}\bigr)
  \longrightarrow \Mcg\bigl(S;B,\{A_i\}\bigr),$$
hence a smooth bundle over $\CP^1_\circ=\CP^1\setminus\{0,1,\omega,\omega^2\}$ with
fiber $S$, one distinguished section corresponding to $B$, and an
unordered three-point multisection corresponding to $\{A_0,A_1,A_2\}$.

  \subsection{{The branched cover of the four-punctured sphere}}

Next we describe a  $288$-fold branched cover of the 2-dimensional sphere
$S$. Fix the four marked points $B, A_0, A_1, A_2\in S$ as before, and let
$S^\circ$ denote their complement $S\setminus \{ B, A_0, A_1, A_2\}$.

In the following we will use the convention that permutations act on the left and are composed right to left:
$(\sigma\rho)(i)=\sigma(\rho(i))$.  A finite group $G$ with a subgroup $H$ acts on
the set $G/H$ of \emph{left} cosets $gH$ by left multiplication,
$ \pi(u)\colon gH\longmapsto (ug)H$ $(u\in G)$,
and with the above composition convention $\pi$ is a homomorphism
$G\to\operatorname{Sym}(G/H)$; the stabilizer of the base coset $H$ is $H$
itself.

Now take the simple group $G=U_3(3)=\PU(3,\mathbb F_9)$, of order $6048$ \cite{Atlas} (see also a short recollection in the appendix),
with its ATLAS standard generators $a,b$, so that
\[
|a|=2,\qquad |b|=6,\qquad |ab|=7 ,
\]
where $|g|$ stands for the order of $g\in G$.
Put $h_7=ab$, of order seven, choose an element $h_3$ of order three with
$h_3h_7h_3^{-1}=h_7^{4}$, and set $H=\langle h_7,h_3\rangle$, of order $21$.
In the appendix we provide one such choice.

Turning back to the branched cover of $S$:
fix the local monodromies at the four marked points
coming from the images of five elements $x,y,p,q,r\in G$ in the symmetric group
$S_{288}=\operatorname{Sym}(G/H)$ with
\begin{equation}\label{eq:groupdata}
    |x|=2,\quad |y|=3,\quad |p|=|q|=|r|=4,
    \qquad
    q=ypy^{-1},\quad r=yqy^{-1},\quad xpqr=1 ,
\end{equation}
and write $X,P,Q,R\in S_{288}$ for the permutations that $x,p,q,r$
induce on $G/H$.  (The elements $x,y,p,q,r$ are given in terms of the
ATLAS generators $a,b$ in the appendix.)  The monodromy sends the
peripheral loops at $B,A_0,A_1,A_2$, taken in that order and
multiplying to $1$, to $X,P,Q,R$.  Their cycle types are
$2^{144},4^{72},4^{72},4^{72}$, so the covering is branched with index
$2$ over $B$ and index $4$ over each $A_i$.

Let $C$ denote the 288-fold branched cover of $S$, branched along the
four points $B, A_0, A_1, A_2$ with monodromy representation given
above by the permutations $X,P,Q,R$ around the four distinguished
points, respectively.  As $x, p, q, r$ generate $G$, the subgroup of
$S_{288}$ generated by these four elements acts transitively on $\{ 1,
\ldots ,288 \}$. Consequently, the branched cover provides a connected
surface, implying that $ \chi(C)=2\cdot288-144-3\cdot216=-216$, so the
surface $C$ has genus $109$.

\begin{rem}\label{rem:torus-fibration}
  As a toy example, a simpler construction can be carried out by taking the 
 \emph{double} branched cover of $S$ branched in the given
 four marked points; the resulting surface is easily seen to be a torus.
 We will follow the constructions with this toy example as well.
  \end{rem}

\subsection{Lifting the monodromy to $C$}\label{ss:lift}

Let $p\colon C\to S$ be the branched covering above and let
$\rho\colon\pi_1(S^\circ,*)\to S_{288}$ have monodromy tuple
$T=(X,P,Q,R)$.  We use a right Hurwitz action:
\[
 T^{\sigma_1}=(X,PQP^{-1},P,R),\qquad
 T^{\sigma_2}=(X,P,QRQ^{-1},Q),\qquad
 T^{\phi\psi}=(T^\phi)^\psi.
\]
Equivalently, $\rho^\phi=\rho\circ\phi_*$ for the corresponding
automorphisms of the fundamental group.  The two positive
half-twists $\sigma _1, \sigma _2$ interchange $(A_0,A_1)$ and $(A_1,A_2)$, respectively.

\begin{lem}\label{lem:conjugator}
  For $i=1,2$ there is a unique $g_i\in S_{288}$ such that
  $T^{\sigma_i}=g_iTg_i^{-1}$.  More explicitly,
  \[
     g_1=\pi(c_1),\qquad g_2=\pi(c_2),\qquad
     c_1=(rx)^4,\qquad c_2=(px)^4,
  \]
  and both $g_1,g_2$ have order three.
\end{lem}

\begin{proof}
  Evaluating the words for $c_1, c_2$  provided in the appendix gives
  \[
  \begin{array}{c|cccc}
       &x&p&q&r\\ \hline
   c_1(\,\cdot\,)c_1^{-1}&x&pqp^{-1}&p&r\\
   c_2(\,\cdot\,)c_2^{-1}&x&p&qrq^{-1}&q
  \end{array}
  \qquad c_1^3=c_2^3=1,\quad c_1,c_2\ne1.
  \]
  Only three entries in each row need separate verification.  For example,
  if $c_1$ fixes $x,r$ and carries $q$ to $p$, then conjugating
  $xpqr=1$ gives $xspr=1$, where $s=c_1pc_1^{-1}$; comparison
  with $xpqr=1$ yields $s=pqp^{-1}$.  The other row is similar.
  Applying the homomorphism $\pi$ proves the required identities.

  The  generators $a,b$ also satisfy
  \[
    a=xq^{-1}p^{-1}r^{-1}x,\qquad b=p^{-1}r^{-1}p^2,
  \]
  and consequently
  \[
    \langle x,p,q,r\rangle=G,\qquad
    \langle X,P,Q,R\rangle=\pi(G).
  \]
  The subgroup $H$ is self-normalizing.  Thus the centralizer of this
  transitive coset action has
  \begin{equation}\label{eq:centralizer}
     C_{\operatorname{Sym}(G/H)}(\pi(G))
          \cong N_G(H)/H=1.
  \end{equation}
  Two simultaneous conjugators differ by an element of this
  centralizer, proving uniqueness.  Finally the nontrivial coset action
  of the simple group $G$ is faithful, so $g_i$ has the same order as
  $c_i$.
\end{proof}

The Hurwitz action commutes with simultaneous conjugation.  Hence every
$\phi\in B_3$ has a unique conjugator $g_\phi$ with
$T^\phi=g_\phi Tg_\phi^{-1}$, and our right-action convention gives
\[
 T^{\phi\psi}
   =g_\phi g_\psi T(g_\phi g_\psi)^{-1},\qquad
 g_{\phi\psi}=g_\phi g_\psi.
\]
This homomorphism is defined on $B_3$.  Its value on the central full
twist is nontrivial, as the next corollary records.

\begin{cor}\label{cor:sigma3}
  The braid $\sigma_1^3$ fixes $T$.
  If $\Delta=\sigma_1\sigma_2\sigma_1$, then
  $g_\Delta^2=X$, so $g_\Delta$ is of order four.
\end{cor}

\begin{proof}
  The first assertion follows from $g_{\sigma_1}^3=1$.
  The full twist $\zeta=\Delta^2=(\sigma_1\sigma_2)^3$ fixes $X$
  and conjugates $P,Q,R$ by their product $PQR=X^{-1}=X$.
  Uniqueness therefore gives $g_\zeta=X$.  Thus
  $g_\Delta^2=X\ne1$ and $g_\Delta^4=1$.
\end{proof}

Choose representatives of the marked mapping classes that are complex
linear in oriented local coordinates near the marked points.  Such
representatives are obtained by
local isotopy; the bundle over $\CP^1_\circ$
can be built with transitions of this form.  The covering-space lifting
criterion applies to these representatives: on $S^\circ$ the
pullback monodromy is conjugate to $\rho$, and the lift extends across
the ramification points smoothly.  Indeed, in the local model $z=w^e$,
a transition $z\mapsto \alpha z$ lifts to $w\mapsto cw$ with $c^e=\alpha$.
The lift is unique because the deck group of
$p\colon C\to S$ is the centralizer in~\eqref{eq:centralizer}.
Lifts of isotopies are unique as well.  In particular, a genuine
order-two representative of $\bar\Delta$ lifts to a diffeomorphism
whose square is a deck transformation and hence is the identity.
This geometric lift and the based sheet conjugator $g_\Delta$ are
different objects: the latter includes the transport used to compare
the chosen sheet labels.

\begin{thm}\label{thm:bundle0}
  The fiberwise coverings assemble into a smooth fiber bundle
  $Z^\circ\to\CP^1_\circ$ with fiber $C$.  It is unique as a
  covering of the given four-pointed sphere bundle.
\end{thm}

\begin{proof}
  Choose local trivializations of the marked sphere bundle.  Each
  transition diffeomorphism lifts by the preceding criterion.  On a
  triple overlap the two compositions cover the same transition
  diffeomorphism and differ by a deck transformation of the whole
  fiber covering, so they agree.  The lifts therefore satisfy the
  cocycle condition.  The same uniqueness identifies any two such
  lifted bundles.
\end{proof}

\subsection{Extension over $1,\omega,\omega^2$}\label{ss:mult}
In the next two subsections we extend the bundle $Z^\circ$ (defined
over $\CP ^1_{\circ}$) to a singular fibration over $\CP ^1$.
Consider first the monodromy $\bar\Delta$ and its conjugates:
these are involutions in the
mapping class group of the four-punctured sphere, which fix $B$ and one of
the $A_i$'s (and interchange the other two).
We start with a local parity observation:

\begin{lem}\label{lem:parity}
Let  $p\colon \widetilde S\longrightarrow S$
be an oriented branched covering of connected surfaces, and suppose
that nontrivial orientation preserving involutions
\[
  g\colon\widetilde S\longrightarrow\widetilde S,
  \qquad
  \iota\colon S\longrightarrow S
\]
satisfy
$  p\circ g=\iota\circ p$.
If \(c\in\widetilde S\) is fixed by $g$, then $a=p(c)$ is fixed by
\(\iota\). If \(p\) has local degree \(e\) at \(c\), then \(e\) is odd.
\end{lem}

\begin{proof}
  It is obvious that $\iota (a)=a$.
Choose small invariant discs around \(c\) and \(a\). On their boundary
circles, both involutions are nontrivial and orientation-preserving,
therefore after parametrizing the circles by
\(\mathbb R/\mathbb Z\), both are given by
\[
  t\longmapsto t+\frac12.
\]

The restriction of \(p\) to the boundary circles is an
orientation-preserving covering of degree \(e\). Let
 $ P\colon\mathbb R\longrightarrow\mathbb R $
be a lift of this circle covering, so that
$  P(t+1)=P(t)+e$.
Equivariance gives
\[
  P\left(t+\frac12\right)
  =P(t)+\frac12+k
\]
for some integer \(k\); the integer is independent of \(t\) by
continuity. Applying this identity twice gives
\[
  P(t+1)=P(t)+1+2k.
\]
Comparing this with \(P(t+1)=P(t)+e\), we obtain
$  e=1+2k$, 
thus $e$ is odd.
\end{proof}

\begin{thm}\label{thm:free}
Let \(\iota\) be the order-two representative of \(\bar\Delta\) that fixes
\(B\) and one of the points \(A_i\), say \(A_1\), and interchanges
\(A_0\) and \(A_2\). Then \(\iota\) lifts uniquely to a diffeomorphism
\[
  g\colon C\longrightarrow C.
\]
The lift satisfies \(g^2=1\), acts freely, and
\(C/\langle g\rangle\) has genus \(55\).
\end{thm}

\begin{proof}
The Hurwitz action of \(\bar\Delta\) carries \(T = (X,P,Q,R) \) to a
simultaneously conjugate tuple. Consequently, the pullback of \(p\) by
\(\iota\) has the same monodromy as \(p\), up to relabeling of the
sheets. The classification of covering spaces over the complement of
the marked points therefore gives a lift $ g\colon C\longrightarrow C$
satisfying $ p\circ g=\iota\circ p$.  Because the peripheral
monodromies have the same cycle types, this lift extends uniquely
across the ramification points.

Any two such lifts differ by a deck transformation of \(p\). The deck
group is trivial, so the lift is unique. Since \(\iota^2=1\),
\[
  p\circ g^2=\iota^2\circ p=p.
\]
Thus \(g^2\) is a deck transformation of \(p\), and therefore
$  g^2=1$.

It remains to prove that \(g\) is free.  A fixed point would lie
above one of the two fixed points $B,A_1$ of $\iota$.
Every point above \(B\) has local degree \(2\), as \(X\) has cycle
type \(2^{144}\). Every point above \(A_1\) has local degree \(4\),
because the corresponding peripheral permutation has cycle type
\(4^{72}\). Thus, in either case, the local degree of \(p\) at \(c\)
would be even,  contradicting Lemma~\ref{lem:parity}. Hence
 $ {\rm {Fix}}(g)=\emptyset$ and the  quotient map
$  C\longrightarrow C/\langle g\rangle $
is an unbranched double covering, the quotient having genus 55 as claimed.
\end{proof}

The above result then implies the following:
\begin{cor}\label{cor:multiple}
  The bundle $Z^\circ$ extends smoothly
  over each of $1,\omega,\omega^2$.  Explicitly, let
  $\Z/2$ act on $C\times D^2$ by $(c,w)\mapsto(g(c),-w)$; the action is free
  because $g$ is free, so the quotient is a smooth $4$-manifold, and
  $(c,w)\mapsto w^2$ descends to a smooth map to $D^2$ whose monodromy
  is $g$ and
  whose central fiber is $C/g$ with multiplicity $2$. \qed
\end{cor}

\begin{rem}
  In our toy example (where the generic fiber is
  the double branched cover of $S$, hence is a torus, see
  Remark~\ref{rem:torus-fibration}) the monodromies around the three points
  $\{ 1, \omega , \omega ^2\}$ are all conjugate to $aba$ (where $a$ and
  $b$ are two Dehn twists along simple closed curves in the torus $T^2$
  intersecting each other in one transverse point). This monodromy is identical
  to the monodromy of the $III$-fiber in Kodaira's classification;
  topologically such a fiber is the union of two spheres with self-intersection
  $(-2)$ intersecting each other tangentially in one point.
  \end{rem}

\subsection{Extension over $0$}\label{ss:cusp}
The extension of the bundle $Z^\circ \to \CP ^1_\circ$ is rather more
complicated over $0$; this will be discussed in the present subsection.
After conjugating the marked fiber, choose a disc $D\subset S$ containing
the two colliding points with peripheral permutations $Q,R$, and put
$F=p^{-1}(D)$ for the covering $p\colon C\to S$.  The local braid is the positive cube of the half
twist interchanging these points.  A lift of the single half twist need
not be supported on $F$: its sheet conjugator is
$\pi((px)^4)$, of order three and cycle type $3^{96}$.  The cube has
trivial sheet conjugator, however, and its lift has a representative
$h_0$ equal to the identity on $C\setminus F$.  We construct its filling
directly from an explicit local branched cover.

\begin{prop}\label{prop:pieces}
  The subgroup $\langle Q,R\rangle$ has order $96$ and acts regularly on
  three orbits of size $96$.  Thus $F$ has three connected components,
  each a degree-$96$ branched cover of $D$.  Each component has genus
  $21$, eight boundary circles, and Euler characteristic $-48$.
  The boundary monodromy is $QR$, of cycle type $12^{24}$, and every
  boundary circle maps to $\partial D$ with degree $12$.
\end{prop}

\begin{proof}
  These orbit and cycle calculations are obtained from the generators
  and words in the appendix: enumerate $G$, form its left cosets of
  $H$, and apply $q,r$ by left multiplication.  On each orbit, both
  $Q$ and $R$ have twenty-four cycles of length four and $QR$ has
  eight cycles of length twelve.  Consequently $
  \chi(F_j)=96-24(4-1)-24(4-1)=-48$. Since $F_j$ has eight boundary
  components, its genus is $21$.
\end{proof}

Using the above result, an explicit description of the extension can be given.

\begin{lem}\label{lem:cuspquotient}
  Put
  \[
    A=\begin{pmatrix}i&0\\0&1\end{pmatrix},\qquad
    B=\frac{1+i}{2}\begin{pmatrix}1&1\\-1&1\end{pmatrix},
  \]
  and
  \[
    F_8(u,v)=u^8+14u^4v^4+v^8,\qquad
    F_{12}(u,v)=u^{12}-33u^8v^4-33u^4v^8+v^{12}.
  \]
  The map
  \[
     \Psi=(F_8,F_{12})\colon\C^2\longrightarrow\C^2
  \]
  is the degree-$96$ quotient by $\Gamma=\langle A,B\rangle$.
  Its branch divisor, in target coordinates $(a,b)$, is $b^2=a^3$,
  with ramification index four along its smooth part.  With suitable
  peripheral paths and sheet labels, its two cusp meridians have the
  same permutation action as $Q,R$ on each orbit of
  Proposition~\ref{prop:pieces}.
\end{lem}

\begin{proof}
  The two polynomials are invariant under $A,B$, and direct calculation
  gives, with $J_6=uv(u^4-v^4)$,
  \begin{equation}\label{eq:twoexpressions}
      F_8^3-F_{12}^2=108J_6^4,  \qquad
      \det\frac{\partial(F_8,F_{12})}{\partial(u,v)}
          =-1728J_6^3.
  \end{equation}
  The projective action of $\Gamma$ on
  $\{0,\infty,1,-1,i,-i\}$ is the octahedral group of order $24$.
  A scalar $\lambda I\in\Gamma$ preserves both polynomials, so
  $\lambda^8=\lambda^{12}=1$ and hence $\lambda^4=1$.
  Conversely, $(AB)^3=iI$.  Thus the scalar kernel is exactly
  $\{1,i,-1,-i\}$ and $|\Gamma|=96$.

  The two invariants have no common nonzero zero.  Indeed, the left
  expression of Equation~\eqref{eq:twoexpressions} would imply
  $J_6=0$; on $u=0$ or $v=0$ one has $F_8\ne0$, and on $u^4=v^4\ne0$
  one has $F_8=16v^8\ne0$.  Hence $\Psi$ is finite of degree
  $8\cdot12=96$.  It factors through $\C^2/\Gamma$, and the induced
  finite map to $\C^2$ has degree one.  Since the target is normal,
  this map is an isomorphism.  The two expressions of
  Equation~\eqref{eq:twoexpressions} identify the branch divisor and
  its ramification index.

  Evaluate $q,r$ from their printed words in the
  $28$-point generators, and enumerate products in $q,r$ and in
  $A,B$ simultaneously.  Both enumerations have $96$ elements, and
  every repeated matrix word gives the same permutation word.
  This verifies an isomorphism sending $q$ to $A$ and $r$ to $B$.

  We also check the geometric meridians.  On $F_8=1$, the path
  \[
    \gamma(x)=(ixv(x),v(x)),\qquad
    v(x)=(x^8+14x^4+1)^{-1/8},\qquad 0\leq x\leq1,
  \]
  maps under $F_{12}$ from $+1$ to $-1$.
By Equation~\eqref{eq:twoexpressions}
  its image is real and has
  no interior critical point, 
  so it parametrizes that interval monotonically.  The endpoints
  $(0,1)$ and $(i/\sqrt2,1/\sqrt2)$ lie on the fixed lines of $A$
  and $B$, respectively.  The transverse eigenvalue of each
  reflection is $i$.  Thus positive meridians based along this
  interval lift to the local transformations $A,B$ on a common
  sheet.  Deck-labeled sheets give the right regular action.
  To use the left regular convention of this paper, compose inversion
  of the labels with the automorphism
  $A\mapsto A^{-1}$, $B\mapsto B^{-1}$; this automorphism is induced
  by the semilinear map $(u,v)\mapsto(\bar u,-\bar v)$.
  The resulting left regular permutations agree with $Q,R$.
\end{proof}

\begin{thm}\label{thm:milnor}
  Each component of $F$, with the local monodromy $h_0$, is the
  Milnor fiber of an ordinary eightfold point, with its positive
  monodromy.  Gluing three copies of this local fibration extends
  $Z^\circ$ smoothly over $0$.
  The central fiber is reduced and has four irreducible components
  whose normalizations have genera $4,4,10,10$ and degrees
  $36,36,108,108$ over the capped downstairs sphere.  It has three
  ordinary eightfold points, with incidence matrix
  \[
     \begin{pmatrix}2&2&2&2\\0&0&4&4\\1&1&3&3\end{pmatrix}.
  \]
  Rows correspond to singular points, columns to irreducible
  components, and entries count local branches.  Its Euler
  characteristic is $-69$.
\end{thm}

\begin{proof}
  Project the cusp family $b^2=a^3$ to the $a$-disc.  Its two branch
  points are $b=\pm a^{3/2}$, so a positive turn around $a=0$ gives
  the positive cube of a half twist.  Lemma~\ref{lem:cuspquotient}
  identifies the prescribed cover on each local orbit with
  $\Psi$, and the lifted projection is simply $a=F_8(u,v)$.
  This polynomial has eight distinct linear factors: the two roots
  of $t^2+14t+1=0$, with $t=(u/v)^4$, are distinct and nonzero.
  Thus its zero fiber is an ordinary eightfold point in a smooth
  surface, and its nearby compact fibers are the corresponding
  Milnor fibers \cite{Milnor}.

  This comparison also fixes the relative-boundary monodromy used
  in gluing.  Choose an outside basepoint on the boundary of the
  $b$-disc.  Its loop around the $a$-disc bounds a disc disjoint
  from the branch curve, so its sheet transport is the identity,
  as for $h_0$.  The coverings therefore agree on the punctured
  neighborhood with this marking, and their collar identifications
  glue to the identity on $C\setminus F$.  There is no need to
  choose a periodic representative fixing the boundary: the
  relative-boundary monodromy is supplied by the Milnor fibration
  itself.

  Topologically the central fiber is obtained by collapsing each
  of the three components of $F$ to a point.  Its four irreducible
  components correspond to the orbits of $\langle X,P\rangle$,
  of sizes $36,36,108,108$.  The branches at the three singular
  points are the $QR$-cycles assigned to their local and
  complementary orbits; this gives the displayed incidence matrix.
  On a complementary orbit of degree $d$, the capped cover has
  branch indices $2,4,12$, and hence
  \[
    \chi(\widetilde C_d)
      =d\left(2-\frac12-\frac34-\frac{11}{12}\right)
      =-\frac d6,\qquad
    g(\widetilde C_d)=1+\frac d{12}.
  \]
  The normalization of the central fiber therefore has Euler
  characteristic $2(-6)+2(-18)=-48$.  Identifying eight preimages
  to one point at each of the three singular points gives
  $\chi(F_0)=-48-3(8-1)=-69$.
 \end{proof}

\begin{rem}
  In the local model, scalar multiplication by $e^{2\pi i/8}$ is a
  periodic representative of the monodromy on the Milnor fiber,
  allowing its boundary circles to rotate.  The representative
  fixing the boundary pointwise differs by a collar adjustment;
  its eighth power is the product of the boundary Dehn twists.
  These are two representatives of the same geometric monodromy
  with different boundary conventions.  The Milnor fillings kill
  every homology class carried by $F$, as is also clear from the contractible local total space.
\end{rem}

\begin{rem}
  In the torus toy model the lifted monodromy is $a^3$, the
  monodromy of the $I_3$-fiber in Kodaira's classification.
\end{rem}

Let $Z$ be the closed four-manifold obtained from $Z^\circ$ by these
four fillings over $\CP^1$.

\begin{cor}\label{cor:euler}
  $\chi(Z)=39$, and there are no further singular fibers.
\end{cor}

\begin{proof}
  The regular fiber has Euler characteristic $-216$, the fiber
  over $0$ has Euler characteristic $-69$, and each reduced
  multiple fiber has Euler characteristic $-108$.  Thus
  \[
     \chi(Z)
       =2(-216)+(-69+216)+3(-108+216)=39.
  \]
  Over $\CP^1_\circ$ the map is a fiber bundle, so no other
  singular fibers occur.
\end{proof}

\subsection{The thirty-six spheres}\label{ss:spheres}

Every braid in the image of $\beta$ fixes the marked point $B$, so for such a
braid $\phi$ the first entry of $T^\phi$ is $X$ itself; since
$T^\phi=g_\phi Tg_\phi^{-1}$, the conjugator $g_\phi$ centralizes $X$ and
therefore permutes the $144$ points of $C$ lying over $B$, which are the
cycles of $X$.

\begin{lem}\label{lem:Q8}
  The sheet conjugator at $z=0$ is trivial.  On the $288$ sheets,
  the three conjugators $d_1,d_2,d_3$ at $1,\omega,\omega^2$ satisfy
  \[
      d_i^2=X\ne1,\qquad d_1d_2d_3=1,
  \]
  and generate the quaternion group $\mathcal Q\cong Q_8$ with
  center $\{1,X\}$.  On the $144$ cycles of $X$ this action factors
  through $\mathcal Q/\{1,X\}\cong(\Z/2)^2$ and is free.  Hence the
  multisection monodromy has $36$ orbits, all of size four.
\end{lem}

\begin{proof}
  Work in $B_3$, with the right Hurwitz action and its homomorphism
  $\phi\mapsto g_\phi$, and put
  \[
      \Delta=\sigma_1\sigma_2\sigma_1,\qquad
      \rho=\sigma_1\sigma_2,\qquad
      \zeta=\Delta^2=\rho^3.
  \]
  The braid representatives of the four monodromies are
  $\sigma_1^3$, $\Delta$, $\rho\Delta\rho^{-1}$ and
  $\rho^2\Delta\rho^{-2}$.  The  identity
  $\Delta\rho=\zeta\sigma_1^{-1}$ gives
  \[
    \sigma_1^3\Delta(\rho\Delta\rho^{-1})
      (\rho^2\Delta\rho^{-2})
       =\sigma_1^3(\Delta\rho)^3\rho^{-3}
       =\zeta^2.
  \]
  The full twist conjugates $P,Q,R$ by $PQR=X^{-1}=X$ and fixes
  $X$, so uniqueness of the sheet conjugator gives $g_\zeta=X$.
  Corollary~\ref{cor:sigma3} gives $g_{\sigma_1^3}=1$.
  Applying the homomorphism to the displayed identity proves
  $d_1d_2d_3=X^2=1$; likewise each $d_i^2=X$.

  The element $X$ commutes with the $d_i$.  Since
  $d_1d_2=d_3^{-1}$, one has $(d_1d_2)^2=X$, and hence
  $d_1d_2=X d_2d_1$.  These are the quaternion relations.
  They give $\mathcal Q\cong Q_8$, because any proper quotient
  of $Q_8$ kills its central involution, whereas $X\ne1$.

  On the cycles of $X$ the center acts trivially.  No $d_i$ can
  preserve one such two-element cycle setwise: its restriction
  to that cycle would have square equal to the identity, whereas
  $d_i^2=X$ interchanges the two elements.  Thus all three
  nonidentity elements of the quotient $(\Z/2)^2$ are
  fixed-point-free, and its orbits have size four.  There are
  $144/4=36$ of them.
\end{proof}

\begin{thm}\label{thm:spheres}
  The $144$ points sweep out a multisection of $Z\to\CP^1$ with $36$ connected
  components $\Sigma_1,\dots,\Sigma_{36}$.  Each is a smoothly embedded sphere,
  they are pairwise disjoint, and each meets every regular fiber $F$
  transversally in four points, all four intersections being positive; in
  particular $\Sigma_i\cdot F=4$.
\end{thm}

\begin{proof}
  The components correspond to the orbits of Lemma \ref{lem:Q8}, so there are
  $36$ of them, each a degree $4$ cover of the base.  It is the regular
  $(\Z/2)^2$-cover branched over $1,\omega,\omega^2$, and at each of those the
  monodromy is a fixed point free involution, hence two points of ramification
  index $2$; so $\chi=4\cdot2-3\cdot2=2$ and the component is a sphere.

  For embeddedness we look at the two local models.  Over the disc at $z=0$ the
  monodromy is trivial, and the $144$ points lie over $B$, which is not one of
  the two colliding marked points, so they avoid the region $F$ of
  Section~\ref{ss:cusp}; the multisection is therefore four disjoint discs there, away
  from the three singular points of the fiber.  Over the disc at $1$, in the
  local model $(C\times D^2)/(\Z/2)$ of Corollary \ref{cor:multiple}, a pair
  $\{p,g(p)\}$ gives the embedded disc $\{p\}\times D^2$, mapping to the base by
  a double branched cover; distinct pairs give disjoint discs because $g$ is
  free.  Distinct components are disjoint because they are disjoint in every
  fiber.  Finally, $\Sigma_i$ meets a regular fiber in the four points of one
  orbit.  Over $\CP^1\setminus\{1,\omega,\omega^2\}$ the projection
  $\Sigma_i\to\CP^1$ is a covering map, so at each of these four points the
  tangent plane of $\Sigma_i$ maps isomorphically to the tangent plane of the
  base and is therefore complementary to the tangent plane of the fiber: the
  intersections are transverse.  Orient $Z$ by the orientations of the fiber
  and of the base, in that order, and orient $\Sigma_i$ by pulling back the
  orientation of the base.  Then at each of the four points the splitting
  $T_pZ=T_pF\oplus T_p\Sigma_i$ is orientation preserving, so every
  intersection counts positively and $\Sigma_i\cdot F=4$.
\end{proof}

  \begin{prop}\label{prop:squareminus1}
    Each sphere $\Sigma_i$ found in the previous theorem
    has self-intersection $-1$, so these are
  exceptional spheres of blow-ups.
    \end{prop}
\begin{proof}
For the map $Z\to \CP^1$ give the base the orbifold structure
\[
   \mathfrak B=S^2(2,2,2),
\]
with order-two points at $1,\omega,\omega^2$; the point $0$ is ordinary.
The downstairs family extends over each order-two point in the form
\[
 (S\times D_w)/\bigl((s,w)\sim(\iota(s),-w)\bigr),
\]
where $\iota$ is the half-turn fixing $B$ and one of the $A_i$.
In an orientation-compatible complex coordinate $u$ centered at $B$ this
identification is $(u,w)\sim(-u,-w)$.
At $0$, use the standard positive cusp-braid model for
$\beta_0=\bar\sigma_1^3$, with the two colliding $A$-points away from $B$.
In particular the sphere family and its distinguished section extend
smoothly near $B$ at $0$. These are auxiliary downstairs local models;
they agree near $B$ with the local fillings already used to construct $Z$.
Let $\mathcal N$ denote the oriented normal line of this distinguished
section, regarded as a complex orbifold line bundle over $\mathfrak B$.

Choose smoothly varying conformal structures on the sphere fibers,
agreeing with the standard structures in these local models.
Such a choice can be made using fiberwise metrics, averaged under the
finite local actions. 
Locally choose an affine uniformizing coordinate $t$
on a sphere with $B$ at infinity, and let $a_1,a_2,a_3$ be the coordinates
of the other three marked points. Put
\[
 D=\prod_{1\le j<k\le3}(a_j-a_k)^2,
 \qquad u=t^{-1},\qquad e=(\partial/\partial u)_B.
\]
The expression $D$ is symmetric in the three points, so their ordering
is not important. After changing affine coordinate to $t'=a t+b$, with $a\ne0$,
one has
\[
 D'=a^6D,\qquad
 u'=\frac{u}{a+bu},\qquad e'=a e.
\]
Consequently $D e^{-6}$ is a well-defined section of $\mathcal N^{-6}$.
It is nowhere zero while the three $A$-points are distinct. In a
uniformizing chart at an order-two point all three are distinct, and the
section extends there without a zero; the isotropy on $\mathcal N^{-6}$ is
trivial. At $0$ we may take the colliding pair to be
$a\pm z^{3/2}$ and the third point to stay away from $a$. Hence, in a local
trivialization,
\[
 D=z^3 h(z),\qquad h(0)\ne0.
\]
Equivalently, the positive braid $\sigma_1^3$ gives discriminant winding
number $+3$. Thus the section of $\mathcal N^{-6}$ has just one zero,
with topological index three. The Euler-number formula for an oriented
line bundle, applied on orbifold charts, gives
\[
 -6\deg_{\mathrm{orb}}\mathcal N
   =\deg_{\mathrm{orb}}(\mathcal N^{-6})=3,
 \qquad
 \deg_{\mathrm{orb}}\mathcal N=-\frac12.
\]

The degree-four map $r_i=\pi|_{\Sigma_i}$ is an unramified orbifold
covering $\Sigma_i\to\mathfrak B$: its two ordinary ramification points
over each cone point have index two, exactly cancelling the order-two
isotropy. Near a point over $B$ the fiberwise branched covering has local
form $u=v^2$. Its leading quadratic term therefore identifies the complex
normal lines by
\[
  \nu_{\Sigma_i/Z}^{\otimes2}\cong r_i^*\mathcal N.
\]
This is an identification coming from the quadratic term, not from the
ordinary differential, which vanishes in the normal direction. It also
holds at the ramification points of $r_i$: in the uniformizing base
coordinate $w$, the local map is
$(v,w)\mapsto(v^2,w)$ into the downstairs orbifold chart above. At $0$
it extends because the section is disjoint from the cusp-filling region.
Taking degrees, with the orientations already fixed in the proof, gives
\[
  2\,\Sigma_i^2
    =\deg(r_i^*\mathcal N)
    =4\deg_{\mathrm{orb}}\mathcal N
    =-2.
\]
It follows that $\Sigma_i^2=-1$ for every $i$.
\end{proof}

\begin{proof}[Proof of Theorem~\ref{thm:A}]
  Theorem~\ref{thm:bundle0} gives the regular genus-$109$ bundle.
  Corollary~\ref{cor:multiple} supplies the three multiplicity-two
  fibers, and Theorem~\ref{thm:milnor} supplies the remaining fiber,
  including its normalization genera, branch incidences, and Euler
  characteristic.  Corollary~\ref{cor:euler} gives $\chi(Z)=39$.
  Finally Theorem~\ref{thm:spheres} and
  Proposition~\ref{prop:squareminus1} give the thirty-six disjoint
  exceptional spheres and their positive degree-four intersections
  with every regular fiber.  Contracting them decreases the Euler
  characteristic by $36$, so $\chi(Z_0)=3$.
\end{proof}

\section{The holomorphic model}\label{sec:B}

We construct the complex surface directly from the pencil
\eqref{eq:pencildata} of the Introduction.

\subsection{The marked pencil}\label{ss:pencil}
Recall that the cubic, the pencil, the point $q_0$ and the line $V$
are given as
\[
   \mathcal C\colon x_0^3=x_1^2x_2,
   \quad
   V\colon x_1=x_2,
   \quad
   q_0=[0:-\tfrac12:1],
   \quad
   L_z\colon 3x_0=z(2x_1+x_2).
\]
The normalization of the cuspidal cubic is parametrized by
$  [x_0:x_1:x_2]=[t^2:t^3:1]$.
The induced map to the pencil parameter and its derivative are
\[
 z(t)=\frac{3t^2}{1+2t^3},\qquad
 z'(t)=\frac{6t(1-t^3)}{(1+2t^3)^2}.
\]
The exceptional parameters are $0,1,\omega,\omega^2$.
The line $L_0$ passes through the cusp $[0:0:1]$, whereas
$L_{\omega^k}$ is tangent to the cubic at
$p_{\omega^k}=[\omega^k:1:1]$.  The line $V$ contains these three
tangency points and misses $q_0$.  At $z=0$ exactly two of the cubic
intersection points approach the cusp; the third tends to $[0:1:0]$,
corresponding to $t=\infty$, where the map $t\mapsto z$ is unramified.
The three poles of $z(t)$ are simple, and the four marked points are
distinct at $z=\infty$.

Blow up $q_0$, and then blow up twice at each tangency point to separate
the cubic from the tangent fiber.  Denote the resulting surface by
$B_7$.  Over each tangency parameter its fiber is a chain
\[
 S_{1,i}+2S_{2,i}+S_{3,i},\qquad
 (S_{1,i}^2,S_{2,i}^2,S_{3,i}^2)=(-2,-1,-2).
\]
Contract the six disjoint curves $S_{1,i},S_{3,i}$, and write
\[
 \mu\colon B_7\longrightarrow R.
\]
The surface $R$ has six rational double points of type $A_1$ and a map
to $\CP^1$ with three fibers of multiplicity two.  Let $\mathcal
C_R,V_R$ denote the images of the strict transforms of the cubic and
$V$.  The distinguished point of a regular fiber is its intersection
with $V_R$; the three further marked points are its intersections with
$\mathcal C_R$.  The exceptional curve over $q_0$ is a different
section and carries no marked point.  For a fixed (finite) parameter
$z$, use the affine fiber coordinate on the fiber over $z$
\[
 [x_0:x_1:x_2]=[z(1+2s/3):1+s:1].
\]
The distinguished point is $s=0$, and the three cubic intersections
satisfy
\begin{equation}\label{eq:cubic}
 8z^3s^3+(36z^3-27)s^2+54(z^3-1)s+27(z^3-1)=0.
\end{equation}

\begin{lem}\label{lem:agreeingmonodromies}
The local marked-sphere monodromy at $0$ is the positive cube of a
half twist.  At each of $1,\omega,\omega^2$ it is an involution fixing
the distinguished point and one cubic intersection, and exchanging
the other two.  After a quadratic base change at a tangency, the
marked family extends with four distinct points.
\end{lem}

\begin{proof}
Near the cusp put $a=x_0$, $b=x_1$, $x_2=1$.  Then
\[
 b^2=a^3,\qquad z=\frac{3a}{1+2b}.
\]
The analytic fiber coordinate $t_{\rm loc}=b/(1+2b)^{3/2}$, with the
root equal to one at the origin, satisfies
$ t_{\rm loc}^2=z^3/27$
on the branch curve.  A positive turn of $z$ therefore turns the two
colliding points through $3\pi$.

At a tangency put $w=z^3$, make the base change $w-1=v^2$, and set
$s=v\xi$.  Dividing \eqref{eq:cubic} by $v^2$ gives
\[
 8(1+v^2)v\xi^3+(9+36v^2)\xi^2+54v\xi+27=0.
\]
In the projective $\xi$-coordinate, the four points at $v=0$ are $B=0$
and $ A=\{i\sqrt3,-i\sqrt3,\infty\}$; they are all distinct.  The
base-change involution is $(\xi,v)\mapsto(-\xi,-v)$, whose central
action has the stated effect.  The map $z\mapsto z^3$ is unramified at
all three tangencies.
\end{proof}

\begin{thm}\label{thm:monodrMatch}
With a suitable marking and distinguished meridian paths, the
monodromies of the marked-sphere family are exactly those of
Lemma~\ref{lem:relation}.
\end{thm}

\begin{proof}
Put $w=z^3$, send $s=0$ to infinity by $t=1/s$, and put $u=t+2/3$.
Equation~\eqref{eq:cubic} becomes
\[
 27(w-1)u^3+9u+2=0.
\]
The double cover branched over the four marked points consequently
has Weierstrass equation
\[
 y^2=u^3+A(w)u+B(w),\qquad
 A(w)=\frac1{3(w-1)},\quad B(w)=\frac2{27(w-1)}.
\]
Its invariant is
\begin{equation}\label{eq:jz}
 j(w)=1728\frac{4A(w)^3}{4A(w)^3+27B(w)^2}
      =\frac{1728}{w}=\frac{1728}{z^3}.
\end{equation}

The moduli orbifold of a sphere with one distinguished point and three
unordered further points has underlying space a sphere with one point
removed, two cone points of orders two and three, and one cusp.  In
the coordinate $w=1728/j$, its cusp is $0$, its point of order two is
$1$, and its point of order three is $\infty$.  Choose peripheral
generators
\[
 T=\bar\sigma_1,\quad D=\bar\Delta,\quad
 R_0=\bar\sigma_1\bar\sigma_2,\qquad
 D^2=R_0^3=1,\quad TDR_0=1.
\]
Equation~\eqref{eq:jz} identifies the moduli map, and
Lemma~\ref{lem:agreeingmonodromies} fixes the local sign and the
order-two stabilizer.  The pullback $w=z^3$ has cusp width three at
zero; at infinity its local degree three removes the order-three
stabilizer.  Choose the meridian paths from the three lifted sectors
in cyclic order.  The corresponding generators of the index-three
subgroup are
\[
 T^3,\qquad D,\qquad R_0DR_0^{-1},\qquad R_0^2DR_0^{-2}.
\]
Their product is $T^3(DR_0)^3=1$.  These are precisely the four
elements in Lemma~\ref{lem:relation}.  Changing the marking or the
paths changes the representation by the corresponding simultaneous
conjugation or Hurwitz move.
\end{proof}

\begin{proof}[Proof of Theorem~\ref{thm:braidmatch}]
The parametrization of the cubic identifies the four exceptional
values and shows that the marked points are otherwise distinct.
Lemma~\ref{lem:agreeingmonodromies} gives the local behavior, and
Theorem~\ref{thm:monodrMatch} identifies the global marked monodromy.
\end{proof}

\subsection{The normalized parameter}\label{ss:cross}

Normalize a regular marked sphere so that its distinguished point
is $0$, two of the other points are $i,-i$, and the third is $\tau$.
The six choices of ordering of the cubic intersections give the six
possible normalizations.

\begin{thm}\label{thm:B}
The parameters satisfy
\[
 z^3=\frac{27\tau^2(\tau^2+1)^2}{(3\tau^2-1)^3}.
\]
This is a map of degree six.  Its fibers over $0,1,\infty$ are,
respectively,
\[
 \{0,i,-i\},\qquad \{i/3,-i/3,\infty\},\qquad
 \{1/\sqrt3,-1/\sqrt3\},
\]
with ramification indices $(2,2,2)$, $(2,2,2)$, and $(3,3)$.
\end{thm}

\begin{proof}
The transformation $t=(s-i)/(2s)$ sends the normalized configuration
to $(\infty,0,1,\lambda)$, where $\lambda=(\tau-i)/(2\tau)$.
The Legendre formula
gives
\[
 j=256\frac{(1-\lambda+\lambda^2)^3}{\lambda^2(1-\lambda)^2}
   =64\frac{(3\tau^2-1)^3}{\tau^2(\tau^2+1)^2}.
\]
Comparing with \eqref{eq:jz} proves the formula.  Its numerator and
denominator are relatively prime and have degree six.  The identity
\[
 27\tau^2(\tau^2+1)^2-(3\tau^2-1)^3=(9\tau^2+1)^2
\]
gives the fiber over one, including its point at infinity.  The
displayed ramification contributes $3+3+4=10=2\cdot6-2$, so there
is no further ramification.
\end{proof}

\subsection{The holomorphic cover and its completion}\label{ss:marked}

We use the following uniqueness principle throughout.  An
automorphism of the base of a connected finite cover has at most
one lift when the deck group is trivial.  A lift exists if its
action preserves the monodromy class.  Relations among such base
automorphisms hold for their lifts.  In our families, lifts over
base trivializations agree on overlaps because their restrictions
to each fiber are unique.  Holomorphic or antiholomorphic
lifts extend to the finite normalization.

\begin{proof}[Proof of Theorem~\ref{thm:holreal}(1)--(2)]
Let $U=\CP^1\setminus\{0,1,\omega,\omega^2\}$.  The tuple
$(X,P,Q,R)$ determines the connected degree-$288$ covers of the
marked fibers over $U$.  Their ramification indices are $2,4,4,4$.
By Theorem~\ref{thm:monodrMatch} and Lemma~\ref{lem:conjugator},
their monodromy class is preserved by transport.  The deck group is
trivial by \eqref{eq:centralizer}, so the local holomorphic covers
glue uniquely to a cover over $U$.  Equivalently, apply Riemann
existence to the finite cover of the branch complement and then
normalize.  Normalizing $R$ in this cover gives a finite normal
compact cover
\[
 f\colon W\longrightarrow R,
 \qquad h\colon W\longrightarrow\CP^1.
\]
We check its smoothness at the exceptional fibers.

At a tangency, the quadratic base change of
Lemma~\ref{lem:agreeingmonodromies} gives four distinct marked
points, hence a smooth family of genus-$109$ curves.  The involution
$(\xi,v)\mapsto(-\xi,-v)$ lifts uniquely.  Its square is the
identity.  On the central curve its only possible fixed points lie
over the two fixed points of $\xi\mapsto-\xi$, where the local
degrees are two and four.  Lemma~\ref{lem:parity} excludes these
fixed points.  The action on the total space is therefore free, and
its quotient is smooth, with a double fiber whose reduction has
genus $55$.  The downstairs quotient has two $A_1$ points and
minimal resolution $S_{1,i}+2S_{2,i}+S_{3,i}$, so this is the
normalization over the corresponding neighborhood in $R$.
The holomorphic family need not be a product; an averaged smooth
connection nevertheless gives the equivariant smooth trivialization
used in Corollary~\ref{cor:multiple}.

At the cusp use the quotient $(a,b)=(F_8,F_{12})$ of
Lemma~\ref{lem:cuspquotient}.  The actual pencil parameter lifts to
\[
 z=\frac{3F_8(u,v)}{1+2F_{12}(u,v)}.
\]
Choose the analytic eighth root
$\lambda=(3/(1+2F_{12}))^{1/8}$ near the origin and put
$(U,V)=(\lambda u,\lambda v)$.  Its derivative at the origin is a
nonzero scalar, and homogeneity gives exactly $z=F_8(U,V)$.  Thus
the normalization is smooth, and the three local orbits give the
three ordinary eightfold points of Theorem~\ref{thm:milnor}.

These models cover every point of $W$.  Over $U$ the covering is
the one used to construct $Z^\circ$.  At the exceptional values we
have identified the same quotient and polynomial fillings, with
their coverings on punctured neighborhoods.  The unique covering
identifications glue to a fiber-preserving diffeomorphism $W\cong Z$.
This proves both assertions.
\end{proof}

\subsection{The exceptional curves and the canonical divisor}\label{ss:holo}

Use the exceptional classes $E_0,E_{i,1},E_{i,2}$ on
$B_7$, where $i=1,2,3$ labels the tangencies and the second index
records the order of the two blowups.  If $H$
denotes the hyperplane class in $\CP ^2$, in this orthogonal basis we have
\begin{align*}
 K_{B_7}&=-3H+E_0+\sum_i(E_{i,1}+E_{i,2}),\\
 \mathcal C'&=3H-\sum_i(E_{i,1}+E_{i,2}),\qquad
 V'=H-\sum_iE_{i,1},\\
 S_{1,i}&=H-E_0-E_{i,1}-E_{i,2},\qquad
 S_{3,i}=E_{i,1}-E_{i,2},\qquad F=H-E_0.
\end{align*}
Here $F$ is the class of a regular fiber, and primes denote strict
transforms.  The rational
pullbacks of the branch divisors are
\begin{equation}\label{eq:branch-pullbacks}
 \mu^*\mathcal C_R=\mathcal C'+\frac12\sum_iS_{1,i},\qquad
 \mu^*V_R=V'+\frac12\sum_iS_{3,i}.
\end{equation}
The coefficients follow by intersecting with each contracted
curve, which has square $-2$.

\begin{prop}\label{prop:Bsquare}
The distinguished section has rational self-intersection
$V_R^2=-\tfrac12$.
\end{prop}

\begin{proof}
The line $V$ meets the cubic transversally at the tangencies and
does not pass through the second infinitely near centers.  Thus
$V'^2=-2$.  It meets each $S_{3,i}$ once and is disjoint from the
$S_{1,i}$.  Equation~\eqref{eq:branch-pullbacks} then gives
$ V_R^2=-2+3(1-1/2)=-1/2$.
\end{proof}

\begin{thm}\label{thm:holo}
Under the diffeomorphism of Theorem~\ref{thm:holreal}(2), the
thirty-six spheres $\Sigma_i$ are disjoint holomorphic rational
curves of self-intersection $-1$.  Their simultaneous contraction
is a smooth compact complex surface $W_0$.
\end{thm}

\begin{proof}
The components of $f^{-1}(V_R)$ are the same thirty-six packets
used in Theorem~\ref{thm:spheres}, hence are the smooth rational
curves $\Sigma_i$.  Since the ramification index along $V_R$ is two,
\[
 f^*V_R=2\sum_{i=1}^{36}\Sigma_i.
\]
Each component maps to $V_R$ with degree four.  The projection
formula and disjointness give
\[
 2\Sigma_i^2=(f^*V_R)\cdot\Sigma_i
            =V_R\cdot f_*\Sigma_i=4V_R^2=-2.
\]
The curves can therefore be contracted holomorphically.  Denote
the contraction by $\pi\colon W\to W_0$.
\end{proof}

\begin{prop}\label{prop:generaltype}
The surfaces $W$ and $W_0$ are projective.  The surface $W_0$ is
minimal of general type, and
\[
 K_{W_0}^2=9,\qquad e(W_0)=3,\qquad
 \chi(\mathcal O_{W_0})=1.
\]
\end{prop}

\begin{proof}
First observe that there are two effective anticanonical divisors
with no common irreducible component:
\[
 -K_{B_7}\sim\sum_iS_{1,i}+2E_0
             \sim2V'+F_{\rm gen}+\sum_iS_{3,i},
\]
where $F_{\rm gen}$ is a generic fiber.  For any irreducible
curve, at least one of these divisors omits it; their intersection
is therefore nonnegative.  Hence $-K_{B_7}$ is nef, and its square
$9-7=2$ makes it big.

For completeness this also proves projectivity of the intermediate
contraction.  The divisor $A=-K_{B_7}+F$ is nef and big, as is
$A-K_{B_7}$.  By the base-point-free theorem
\cite[Theorem~3.3]{KM}, a multiple of $A$ defines a projective
contraction.  Its contracted curves are exactly the six
$S_{1,i},S_{3,i}$: a horizontal curve has positive intersection
with $F$, while among the fiber components these six are precisely
the ones of anticanonical degree zero.  The contraction is
therefore $R$.  Its finite cover $W$ is projective, as is the
blow-down $W_0$.

The ramification formula in codimension one gives
\[
 K_W=f^*D,\qquad D=K_R+\frac34\mathcal C_R+\frac12V_R.
\]
There is no further divisorial ramification.  Using
\eqref{eq:branch-pullbacks}  we obtain
\[
 \mu^*D=\frac18\left(7H-E_0-3\sum_i(E_{i,1}+E_{i,2})\right),
 \qquad D^2=-\frac3{32},\qquad K_W^2=-27.
\]
Since $K_W=\pi^*K_{W_0}+\sum_i\Sigma_i$, the exceptional divisor
subtracts precisely $\frac12f^*V_R$.  Consequently
\[
 \pi^*K_{W_0}=f^*D_0,\qquad
 D_0=K_R+\frac34\mathcal C_R,\qquad
 \mu^*D_0=-\frac18K_{B_7}.
\]
Thus $K_{W_0}$ is nef and big.  It follows that $W_0$ is of
general type; nefness excludes any further $(-1)$-curve.  Moreover
\[
 K_{W_0}^2=288\frac{2}{64}=9,
 \qquad e(W_0)=39-36=3
\]
by Corollary~\ref{cor:euler}.  Noether's formula gives
$\chi(\mathcal O_{W_0})=(9+3)/12=1$.
\end{proof}

\subsection{The real involution}\label{ss:anti}

We collect the finite data needed for the real structure.  These
statements refer to the printed left coset action, independently
of an identification of $W_0$ with any previously known surface.

\begin{lem}\label{lem:finite-covering-data}
\begin{enumerate}
\item There is a unique permutation $\varsigma$ satisfying
\[
 \varsigma X\varsigma^{-1}=X^{-1},\quad
 \varsigma P\varsigma^{-1}=R^{-1},\quad
 \varsigma Q\varsigma^{-1}=Q^{-1},\quad
 \varsigma R\varsigma^{-1}=P^{-1}.
\]
It has cycle type $1^8 2^{140}$, and
\[
 \varsigma\pi(c_1)\varsigma^{-1}=\pi(c_2)^{-1},\qquad
 \varsigma\pi(c_2)\varsigma^{-1}=\pi(c_1)^{-1}.
\]
\item On the thirty-six four-element packets of $X$-cycles,
$\varsigma$ has type $1^4 2^{16}$.  Each fixed packet contains two
fixed $X$-cycles and two fixed sheets.
\item Put $Q'=PQP^{-1}$.  The tuple $(X,Q',P,R)$ is a Hurwitz
transform and satisfies $XQ'PR=1$.  The three
$\langle P,R\rangle$-orbits are invariant under $\varsigma$.
Its action on the twenty-four $PR$-cycles has type $1^8 2^8$,
with $0,4,4$ fixed cycles on the three local orbits.
\end{enumerate}
\end{lem}
\begin{proof}
  Consider $\varsigma$ given in Equation~\eqref{eq:varsigma} in
  the appendix. Direct computation checks the properties of
  this choice of $\varsigma$ listed in the statement.
Uniqueness follows directly
from \eqref{eq:centralizer}; the square of $\varsigma$ centralizes
the tuple and is therefore the identity.  The identities with
$\pi(c_1),\pi(c_2)$ check compatibility with transport.  Finally,
$PR=(XQ')^{-1}$ places its cycles in the complementary
$\langle X,Q'\rangle$-orbits as well, so the last count is
compatible with the global central fiber.
\end{proof}

\begin{lem}\label{lem:realmodel}
Complex conjugation induces an antiholomorphic involution $\kappa$
of $R$ covering $z\mapsto\bar z$ and preserving the branch divisors.
It preserves the fibers at zero and one and exchanges those at
$\omega,\omega^2$.
\end{lem}

\begin{proof}
The equations of the cubic, line and pencil are real.  The centers
$q_0,p_1$ and the tangent direction at $p_1$ are real; the other
centers and their infinitely near centers form conjugate pairs.
Conjugation therefore lifts to $B_7$, permutes its six contracted
curves, and descends to $R$.
\end{proof}

\begin{thm}\label{thm:anti}
The involution $\kappa$ lifts uniquely to an antiholomorphic
involution $\tilde\kappa$ of $W$.  It permutes the exceptional
curves and descends to a real structure $c$ on $W_0$.  Exactly four
of the exceptional curves have real points, and
\[
 \chi\bigl(W^{\tilde\kappa}\bigr)=-5,\qquad
 \chi(W_0^c)=-1.
\]
\end{thm}

\begin{proof}
Choose the real normalization
\[
 (B,P,Q,R)=(0,i,\tau,-i),\qquad
 \tau\in\R,\quad 0<|\tau|\ll1,
\]
where the letters also indicate the corresponding peripheral
permutations.  Conjugation reverses the peripheral orientations
and exchanges the points at $i,-i$.  Its action on the tuple is
therefore $(X^{-1},R^{-1},Q^{-1},P^{-1})$.  The first part of
Lemma~\ref{lem:finite-covering-data} gives its lift and its
compatibility with base transport.  Fiberwise uniqueness glues
these lifts over the regular base, and normality extends the lift
antiholomorphically over the finite cover.  Its square is the
identity by uniqueness.  Since it preserves $f^{-1}(V_R)$, it
descends through the contraction to $c$.

We compute the real Euler characteristic directly from $h$.
Theorem~\ref{thm:B} gives $z^3<0$ in the chosen normalization.
As $\tau\to0$, the points $B,Q$ coalesce in this coordinate, and
the complementary vanishing region of the original cubic cusp
contains the points with meridians $P,R$.  The Hurwitz move to
$(X,PQP^{-1},P,R)$ places those meridians together.  Its three
local orbits are invariant under $\varsigma$, and their boundary
cycles have $0,4,4$ fixed elements.  Thus all three singular points
are real, with eight individually real branches in total.

The real locus of the normalization of the central curve is a
disjoint union of circles, possibly empty.  Gluing the branches
removes eight real normalization points and adds three real
singular points.  Hence
\[
 \chi\bigl((h^{-1}(0))^{\tilde\kappa}\bigr)=3-8=-5.
\]
Every other real fiber has real locus a union of circles, including
the reduction of the double fiber over one.  Additivity of Euler
characteristic for the proper real-algebraic map to $\R P^1$ now
gives $\chi(W^{\tilde\kappa})=-5$.

The packet action has four fixed packets, each containing fixed
points.  The corresponding invariant rational curves therefore
have real circles.  Contracting each circle increases the real
Euler characteristic by one; the remaining sixteen conjugate
pairs do not change the real locus.  Thus
$\chi(W_0^c)=-5+4=-1$.
\end{proof}

\begin{proof}[Completion of the proof of Theorem~\ref{thm:holreal}]
Theorem~\ref{thm:anti} proves assertion (3), completing the proof.
\end{proof}

\begin{proof}[Proof of Theorem~\ref{thm:CSid}]
Proposition~\ref{prop:generaltype} proves that $W_0$ is minimal of
general type with $K^2=9=3e$.  Equality in the Miyaoka--Yau
inequality makes it a smooth ball quotient \cite{Mi,Ya}.  The
classification of surfaces of general type of Euler number three
leaves the fake projective planes and the Cartwright--Steger
surface \cite{PY,CS,BY}.  Fake projective planes admit no real
structure \cite{KK}; Theorem~\ref{thm:anti} therefore identifies
$W_0$ with the Cartwright--Steger surface $Y$.

Consequently $W$ is its blow-up at the thirty-six images of the
curves $\Sigma_i$.  We henceforth write $\CSt=W$ and $\CS=W_0$
under this identification.  It transports $c$ to a real structure
on $Y$.  To compare it with the Borisov--Yeung involution $c_0$,
use $\operatorname{Aut}(Y)\cong C_3$ \cite[Main Theorem (e)]{CKY}.
Any other real involution is $a c_0$ for $a\in C_3$.
Conjugation by $c_0$ acts on $C_3$ either trivially or by inversion.
In the first case $(a c_0)^2=1$ forces $a=1$.  In the second,
choose $b\in C_3$ with $b^2=a$; then
$b c_0 b^{-1}=b^2c_0=a c_0$.  Thus the real structures agree up
to a holomorphic automorphism of $Y$.
\end{proof}

\begin{lem}\label{lem:ordthree}
The automorphism $[x_0:x_1:x_2]\mapsto[\omega x_0:x_1:x_2]$
lifts uniquely to a holomorphic automorphism $\Phi$ of $W$ covering
$z\mapsto\omega z$.  It has order three and descends to a generator
$\sigma$ of $\operatorname{Aut}(Y)$.  Moreover
\[
 \tilde\kappa\Phi\tilde\kappa=\Phi^{-1},\qquad
 c\sigma c=\sigma^{-1}.
\]
\end{lem}

\begin{proof}
The projective automorphism preserves the pencil, the branch
divisors and all blowup centers as sets.  Its cyclic motion of the
three cubic intersections preserves the covering class by
Lemma~\ref{lem:conjugator}.  The same fiberwise uniqueness and
normalization argument used for conjugation gives the holomorphic
lift.  Its cube covers the identity of $R$, so it is a deck
transformation and is the identity.  The lift is nontrivial because
it moves the base.  It permutes the exceptional curves and hence
descends to a nontrivial order-three automorphism of $Y$; by
\cite[Main Theorem (e)]{CKY}, this generates its automorphism group.
Conjugation takes the downstairs multiplier $\omega$ to
$\omega^{-1}$, and uniqueness of lifting gives both displayed
relations.
\end{proof}

\begin{rem}\label{rem:orderthreeorder}
The sheet conjugator for the braid $\varrho=\sigma_1\sigma_2$
has order six, since $g_{\varrho^3}=X$.  This is consistent with
$\Phi^3=1$: a based braid includes the transport used to compare
two fibers, whereas $\Phi$ is an automorphism of the entire
covering family.
\end{rem}

Note that the identification of $W_0$ (and so $Z_0$) with
the Cartwright-Steger surface implies that
\[
H_1(W;\Z)\cong H_1(W_0; \Z)\cong H_1(Z; \Z )\cong H_1 (Z_0; \Z )\cong\Z^2.
\]

\appendix

\section{Group theoretic properties of $U_3(3)$}

\newcommand{\PSU}{\operatorname{PSU}}
\newcommand{\PGU}{\operatorname{PGU}}
\newcommand{\SU}{\operatorname{SU}}
\newcommand{\GU}{\operatorname{GU}}
\newcommand{\Aut}{\operatorname{Aut}}
\newcommand{\Out}{\operatorname{Out}}

\subsection*{The simple unitary group $U_3(3)$}
Let $\F_9/\F_3$ be the quadratic extension of $\F_3$ and let
$\sigma \colon x\mapsto x^3$ be its nontrivial Galois automorphism.  On
the vector space $V=\F_9^3$,
fix a nondegenerate Hermitian form, for example
\[
 h(v,w)=v_1w_1^\sigma+v_2w_2^\sigma+v_3w_3^\sigma.
\]
The full and special unitary groups are
\[
 \GU_3(3)=\{g\in \operatorname{GL}_3(9)\mid h(gv,gw)=h(v,w)\},
 \qquad
 \SU_3(3)=\GU_3(3)\cap \operatorname{SL}_3(9).
\]
Their projective images are denoted by $\PGU_3(3)$ and $\PSU_3(3)$.
For general $q$ one has
\[
 |\SU_3(q)|=q^3(q^2-1)(q^3+1),
 \qquad
 |Z(\SU_3(q))|=\gcd(3,q+1).
\]
Consequently,
\[
 |\PSU_3(3)|=3^3(3^2-1)(3^3+1)=27\cdot 8\cdot 28=6048.
\]
Since $\gcd(3,4)=1$, the center of $\SU_3(3)$ is trivial.  Moreover,
$\PGU_3(3)/\PSU_3(3)$ has order $\gcd(3,4)=1$.  Thus in this particular case the natural projective map gives
\[
 \SU_3(3)\xrightarrow{\sim}\PSU_3(3)=\PGU_3(3).
\]
After making this natural identification, the common abstract group is often
written $U_3(3)$ or $\mathrm{PU}(3,\F_9)$.
The general simplicity theorem for projective special unitary groups implies
that $\PSU_3(3)$ is nonabelian simple; the exceptional nonsimple case in
rank $2$ is $\PSU_3(2)$ \cite{carter,taylor}.

The isotropic one-dimensional subspaces of $V=\F_9^3$ form the Hermitian curve
\[
 \mathcal H(2,9)=
 \bigl\{[x:y:z]\in \operatorname{PG}(2,9)\mid x^4+y^4+z^4=0\bigr\}.
\]
For a Hermitian curve over $\F_{q^2}$ the number of isotropic points is
$q^3+1$; hence $|\mathcal H(2,9)|=28$.  The group $G=U_3(3)$ acts
faithfully and doubly transitively on these points.  A point stabilizer has
order
$ |G|/28=216.$

\subsection*{Generators and presentations of $U_3(3)$}
The ATLAS standard generators are elements $a,b\in G$ satisfying
\[
 \ord(a)=2,\qquad \ord(b)=6,\qquad \ord(ab)=7.
\]
Thus $G$ is a quotient of the hyperbolic triangle group
\[
 \Delta(2,6,7)=\langle x,y\mid x^2=y^6=(xy)^7=1\rangle,
\]
where hyperbolicity follows from
$\frac12+\frac16+\frac17<1$.  The three order conditions alone do not
characterize the finite quotient.  A presentation on the standard pair is
\[
\begin{split}
G\cong \langle a,b\mid {}&a^2=b^6=(ab)^7=1,\\
 &[a,(ab^2)^3]=1,\qquad
 b^3[b^2,ab^3a]^2=1\rangle,
\end{split}
\]
where $[x,y]=x^{-1}y^{-1}xy$ \cite{onlineatlas}.

A particularly useful concrete realization is the following primitive
permutation representation on $28$ points:
\[
\begin{aligned}
a={}&(2\,3)(4\,6)(5\,8)(7\,11)(9\,13)(10\,15)\\
   &\quad (12\,14)(16\,20)(17\,22)(18\,23)(24\,27)(25\,28),\\[1mm]
b={}&(1\,2\,4\,7\,12\,17)(3\,5\,9\,14\,19\,22)\\
   &\quad (6\,10\,13\,18\,24\,23)(8\,11\,16\,21\,26\,28)
          (20\,25\,27).
\end{aligned}
\]
Direct calculation gives
\[
 |\langle a,b\rangle|=6048,
 \qquad
 \bigl(\ord(a),\ord(b),\ord(ab)\bigr)=(2,6,7).
\]
Since $G$ is simple, every nontrivial permutation action is faithful.
The displayed degree-$28$ action is a convenient faithful model for 
finite calculations.

The elements used in the representations defining the branched covers
can be explicitly described in the
ATLAS presentation of $G$:
writing $B=b^{-1}$ and using $a^{-1}=a$, the five elements of
Equation~\eqref{eq:groupdata} are
\[
 \begin{aligned}
   h_3&=aBababbaBaBaB, &\qquad y&=abaB, &\qquad x&=abaBBaBaBabababa,\\
   p&=aBaBaBabababaBaB, &\qquad q&=bbbabbaBabbabb, &\qquad
   r&=aBabababbababa ,
 \end{aligned}
\]
and $X,P,Q,R$ are their images under $\pi\colon G\to\operatorname{Sym}(G/H)$.
The element $\varsigma$ of Lemma~\ref{lem:finite-covering-data}
can be given in this presentation as follows.
Define the following permutation of the 28-point set:
\[
\theta = (1\, 22)(2\, 19)(3\, 17)(4\, 14)(5\, 12)(6\, 18)
(7\, 9)(8\, 16)(10\, 13)(20\, 27)(21\, 28)(23\, 24).
\]
The action $\varsigma$ on a coset $gH$ is then given
by
\begin{equation}\label{eq:varsigma}
  \varsigma (gH) =  \theta g\theta ^{-1}baH.
  \end{equation}
Direct calculation then verifies the claims of
Lemma~\ref{lem:finite-covering-data}.

\newcommand{\w}{w}
\subsection*{Matrix representations}
These elements can then be represented by $3\times 3$ matrices over the field
of 9 elements as follows. Let
\[
 \F_9=\F_3[\w]/(\w^2-\w-1).
\]
Thus $\w^2=\w+1$, the element $\w$ has multiplicative order $8$, and the
nontrivial automorphism of $\F_9/\F_3$ is
$ \sigma\colon x\longmapsto x^3$.

Let the matrices $A,B$ correspond to the generators $a,b$ of the group
$G=U_3(3)$; with this notation
\[
 A=
 \begin{pmatrix}
  1&\w^{5}&\w^{2}\\
  \w^{7}&0&\w^{5}\\
  \w^{6}&\w^{7}&1
 \end{pmatrix},
 \qquad
 B=
 \begin{pmatrix}
  \w^{5}&\w^{4}&1\\
  \w^{3}&\w^{6}&\w^{2}\\
  0&\w^{5}&\w^{5}
 \end{pmatrix}.
\]
These are the standard generators in the natural, absolutely irreducible
$3$-dimensional ATLAS representation over $\F_9$.
A direct calculation gives
\[
 \ord(A)=2,
 \qquad \ord(B)=6,
 \qquad \ord(AB)=7.
\]
Moreover, direct enumeration gives
\[
 |\langle A,B\rangle|=6048.
\]
Since $\langle A,B\rangle\leq U_3(3)$ and $|U_3(3)|=6048$, it follows
that $ \langle A,B\rangle= U_3(3)$.  For later reference, the product
of order $7$ is
\[
 AB=
 \begin{pmatrix}
  \w^{3}&\w^{4}&\w^{6}\\
  \w^{4}&\w^{4}&\w^{5}\\
  \w^{4}&\w^{3}&\w^{6}
 \end{pmatrix}.
\]

\end{document}